\input amstex
\input amsppt.sty
\input epsf
\magnification=1200
\parindent 20 pt
\def\nologo{\let\log@\empty}
\nologo
\vsize=7.50in
\NoBlackBoxes
 \define\pf{\demo{Proof}}
\define\tZ{\tilde Z}\define\tY{\tilde Y}
\define\tX{\tilde X}
\define\BZ{\Bbb Z}
\define\Bp{\Bbb P}
\define\Center{\operatorname{Center}}
\define\tB{\tilde B}
\define\tT{\tilde T}

\define\ep{\endproclaim}
\define \a{\alpha}
\define \be{\beta}

\define \g{\gamma}

\define \s{\sigma}
\define \fa{\forall}

\define \la{\langle}
\define \ra{\rangle}

\define \G{\Gamma}

\define \CPt{\Bbb C\Bbb P^2}

\define \BR{\Bbb R}

\define \vp{\varphi}

\define \Dl{\Delta}
\define \dl{\delta}
\define \BC{\Bbb C}

\define \1{^{-1}}\define \2{^{-2}}
\define \p{\partial}

\define \Int{\operatorname{Int}}
\define \un{\underline}
\define \ov{\overline}

\define \df{\dsize\frac}
\define \fc{\frac}

\define \edm{\enddemo}

\define \bk{\bigskip}

\topmatter

\title The Fundamental Group  of a  $\CPt$  Complement of a Branch Curve
 \\
  as an Extension of a Solvable Group by a
 Symmetric Group\endtitle
\author Mina Teicher\endauthor
\affil Department of Mathematics, Bar-Ilan University, 52900 Ramat-Gan,
Israel\endaffil
\email teicher\@macs.biu.ac.il\endemail
\leftheadtext{M. Teicher}
\rightheadtext{Fundamental group related to Veronese branch curve}
\thanks  This research was partially supported by the Emmy Noether Research
Institute
 of Bar-Ilan University and the Minerva Foundation of Germany.\endthanks
\abstract
The main result in this paper is as follows:
\proclaim{Theorem} Let $S$ be the branch curve in $\Bbb C\Bbb P^2$ of a generic
projection of a Veronese surface.  Then $\pi_1(\Bbb C\Bbb P^2-S)$ is an
extension of a solvable group by a symmetric group.\endproclaim

A group with the property mentioned in the theorem is ``almost solvable'' in
the sense that it contains a solvable normal subgroup of finite index.
We pose the following question.
\proclaim{Question}
For which families of simply connected  algebraic surfaces of
general type is the fundamental group of the complement of the branch curve of
a generic projection to  $\Bbb C\Bbb P^2$  an extension of a solvable group by
a symmetric group?\endproclaim
\endabstract
\endtopmatter

\document

\subheading{Introduction}
\baselineskip 19pt
Our study of fundamental groups of complements of branch curves is part of our
plan to use fundamental groups in order to distinguish among different
components of moduli spaces of surfaces of general type.

There are not many known computations of fundamental groups of complements of
branch curves.
The topic started with Zariski who proved
in the 30's that if  $X$ is a cubic surface in $\Bbb C\Bp^3$ and $S$ is the
branch curve of a generic projection of $X$ then
$\pi_1(\CPt-S)\simeq Z_2\star Z_3$ (see \cite{Z}).
In the late 70's Moishezon proved that if $X$ is a $\deg n$ surface in
$\Bbb C\Bp^3$ then  $\pi_1(\CPt-S)\simeq B_n/\Center,$ where  $B_n$ is the
braid
group  of order $n $ (see \cite{Mo}.
In fact, Moishezon's result for
$n=3$ is the same as Zariski's result since  $B_3/\Center\simeq Z_2\star Z_3.$
The next example was Veronese of order 2 (see \cite{MoTe3}).
In all the above examples   $\pi_1(\CPt-S)$ contains a free
noncommutative subgroup with 2 generators, so it is ``big''.

Unlike in the early results, in  this   paper we present
    $\pi_1(\CPt-S)$
not ``big''.
We study here the fundamental group of the complement in
$\CPt$ of the branch curve of a generic projection of a Veronese surface.

Our main result is as follows:
\proclaim{Theorem} Let $S$ be the branch curve in $\CPt$ of a generic
projection of a Veronese surface.  Then $\pi_1(\CPt-S)$ is an extension of a
solvable group by a symmetric group.\ep

We believe that the statement of the theorem is valid for many classes of
surfaces of general type.
A group with the property mentioned in the theorem is ``almost solvable'' in
the sense that it contains a solvable normal subgroup of finite index.
We pose the following question.
\proclaim{Question}
For which families of simply connected  algebraic surfaces of
general type is the fundamental group of the complement of the branch curve of
a generic projection to  $\CPt$  an extension of a solvable group by a
symmetric group?\ep

In \cite{MoTe10}, Proposition 2.4, we proved an almost solvability
theorem for the complement of $ S$ in $\BC^2 $ where $\BC^2$ is a generic
affine piece of
$\CPt.$
In this paper we move from $\BC^2$ to $\CPt$.
This situation
involves new techniques (\S3 and \S5),  the Van Kampen Theorem for projective
curves, quoted in
\S1, and different results on the structure of
$\pi_1(\BC^2-S)$ from \cite{MoTe9} and \cite{MoTe10} quoted in \S4.
To formulate the results in \S3 and \S4 we need some information on the braid
group
$B_n$ and its quotient $\tB_n$  which we give in \S2.
The main theorem is proven in \S5.

The theorem can be generalized for any Veronese embedding.
In this paper we choose to prove it for an embedding of order 3 to simplify the
presentation.
The result for any Veronese appears in Section 6.

The  braid group $B_n$ plays an important role in describing fundamental groups
of complements of curves.
There is a quotient of $B_n,$ namely $\tB_n$, which acts on our group
$\pi_1(\CPt-S)$.
We believe that $\tB_n$ acts on fundamental groups of complements of branch
curves for many classes of surfaces of general type, and we can characterize
such fundamental groups through the classification of $\tB_n$-groups.

Lately, there is a growing interest in fundamental groups in general, in
classical algebraic geometry  and in K\"ahler geometry cf., for
example,
\cite{L}, \cite{Si}, \cite{To}.  For fundamental groups of complements of
curves see also \cite{CT} and \cite{DOZ}.
\bk

\subheading{\S1. The Van Kampen Theorem}

As we stated in the introduction, our starting point for proving the main
theorem is the Van Kampen Theorem for projective complements of curves.
The Van Kampen Theorem from the 1930's deals with fundamental groups of affine
and projective complements of curves.
Since our main result is a statement on the fundamental group of the complement
in
$\CPt$, we shall  only quote in this section the  Van Kampen Theorem for the
projective complement.
We shall start with a few definitions, that we need in order to formulate
the theorem.

\definition{1.1 Definition}\ $\underline{\ell(\g)}$.

Let $D$ be a disk.
Let $p\in \Int(D).$
Let $u\in\partial D.$
 Let $\g$ be a simple path connecting $u$
with $p$.
We assign to $\g$ a loop   as follows:
Let $c$ be a small (oriented) circle around $p.$
Let $\g'$ be the part of $\g$ outside of $c.$
We define
 $\ell(\g)=\g^{\prime}\cup
c\cup
\g^{\prime\,-1}.$ We also use the same notation $\ell(\g)$   for the element of
$\pi_1(D-K,u)$ corresponding to $\ell(\g)$.
If $p\in K,$\ $K\subset D,$\ $K$ finite, and $\g$ does not meet any
other point of $K,$ then $\ell(\g)$ can be chosen to be in
$\pi_1(D-K,u).$\enddefinition

\definition{1.2 Definition} \ $\underline{g\text{-base (good geometric base)}}$

Let $D$ be a disk, $K\subseteq D,$ $ K=\{a_1,\dots,a_m\}.$
Let $u\in \p (D)-K.$
Let $\{\g_i\}_{i=1}^m$ be a bush in $(D,K,u),$ i.e., $\g _i$ is a simple path
connecting $u$ with $a_i$,\
$\forall i,j$\
$\g_i\cap\g_j=u,$ \ $\forall i$\ $\g_i\cap K$ = one point, and
$\{\g_i\}$ are ordered counterclockwise around $u.$
Let $\G_i=\ell (\g_i) \in \pi_1 (D-K,u)$ be the loop around $a_i$
determined by $\g_i.$\
$\{\G_i\}_{i=1}^m$ is a $g$-base of $\pi_1(D-K, u).$\enddefinition

\definition{1.3 Remark} A $g$-base is a free base of $\pi_1(D-K,*)$ which
is essential
in the formulation of the  Van Kampen Theorem.\enddefinition

\subheading{1.4}\ Consider the following situation:  Let $S$ be a curve in
$\CPt$ of $\deg m$, s.t. $S$ is transversal to the line in infinity.
Let $\pi: \BC^2\to\BC$ be a generic projection.
Let $N=\{x\in \BC\bigm| \#\pi\1(x)\cap S\underset \neq\to < m\}.$
Let $u\in \BC-N$, s.t. $u$ is real and $|x|<u$\ $\forall x\in N.$
Let $\BC_u=\pi\1(u).$
Let $\{\G_i\}_{i=1}^m$ be a $g$-base of $\pi_1(\BC_u-S\cap \BC_u,*)$.
By abuse of notation we also use the notation $\G_i$ for the image of $\G_i$
in $\pi_1(\BC^2-S,*).$

\proclaim{1.5 Theorem} \ \rom{(Projective Van Kampen Theorem)}
\quad In the situation of $1.4$ we have
$$\pi_1(\CPt-S,*)\simeq
\fc{\pi(\BC^2-S,*)}{\big\langle\prod\limits_{i=1}^m\G_i\big\rangle}.$$
where $\big\langle\prod\limits_{i=1}^m \G_i \big\rangle$ is the subgroup
normally
generated by $ \prod\limits_{i=1}^m\G_i .$\ep

\demo{Proof} \ \cite{VK}\edm

\bk
\subheading{\S2.\ Introducing $\tB_n$, a quotient of $B_n$}

In this section we bring the definition of the braid group and we distinguish
certain elements, called half-twists.
Using half-twists we present Artin's Structure Theorem for the braid group and
the natural homomorphism to the symmetric group.
We also define transversal half-twists and the quotient of $B_n$ called
$\tB_n.$

\definition{2.1 Definition} $\un{\text{Braid group}\ B_n=B_n[D,K]}$

Let $D$ be
a closed disc in $\Bbb R^2,$ \ $K\subset D,$ \ $K$ finite.
Let $B$ be the group of
all diffeomorphisms $\beta$ of $D$ such that $\beta(K) = K\,,\, \beta
|_{\partial D} = \text{Id}_{\partial D}$\,.
For $\beta_1 ,\beta_2\in B$\,, we
say that $\beta_1$ is equivalent to $\beta_2$ if $\beta_1$ and $\beta_2$ induce
the same automorphism of $\pi_1(D-K,u)$\,.
The quotient of $B$ by this
equivalence relation is called the braid group $B_n[D,K]$ ($n= \#K$).
The elements of $B_n[D,K]$ are called braids.
\enddefinition

\newpage

\definition{2.2 Definition}\ \underbar{$H(\sigma)$, half-twist defined by
$\sigma$}

Let
$D,K$ be as above.
Let $a,b\in K\,,\, K_{a,b}=K-a-b$ and $\sigma$ be a simple
path in $D-\partial D$ connecting $a$ with $b$ s.t. $\sigma\cap
K=\{a,b\}.$
Choose a small regular neighborhood $U$ of $\sigma$ and an
orientation preserving diffeomorphism $f:{\Bbb R}^2 \rightarrow {\Bbb C}^1$\
(${\Bbb C}^1$ is taken with the usual ``complex'' orientation) such that
$f(\sigma)=[-1,1]\,,\,$ \ $ f(U)=\{z\in{\Bbb C}^1 \,|\,|z|<2\}$\,.
Let $\alpha(r),r\ge 0$\,,
be a real smooth monotone function such that $
\alpha(r) = 1$ for $r\in [0,\tsize{3\over 2}]$ and
                $\alpha(r) =   0$ for $ r\ge 2.$

Define a diffeomorphism $h:{\Bbb C}^1 \longrightarrow {\Bbb C}^1$ as follows.
For $z\in {\Bbb C}^1\,,\, z= re^{i\varphi},$ let
$h(z) =re^{i(\varphi +\alpha(r))}$\,.
It is clear that on
$\{z\in{\Bbb C}^1\,|\,|z|\leq\tsize{3\over 2}\}$,\ $h(z)$ is the
positive rotation by $180^{\tsize{\circ}}$ and that
$h(z)=\text{Identity on }\{z\in{\Bbb C}^1\,|\,|z|\ge 2\}$\,, in
particular, on ${\Bbb C}^1 -f(U)$\,.
Considering $(f\circ h\circ f^{-1})|_{D}$ (we always take composition from
left to right), we get a diffeomorphism of $D$ which interchanges $a$ and
$b$ and
is the identity on $D-U$\,.
Thus it defines an element of $B_n[D,K],$ called the
half-twist defined by $\sigma$ and denoted $H(\sigma).$ \enddefinition

Using half-twists we build a set of generators for $B_n.$
 \definition{2.3 Definition}\  $\underline{\text{Frame of}\ B_n[D,K]}$

Let $D$ be a disc in $\BR^2.$
Let $K=\{a_1,\ldots ,a_n\},$\ $K\subset D.$
 Let $\sigma_1,\ldots ,\sigma_{n-1}$
be a system of simple paths in $D-\partial D$ such that each $\sigma_i$
connects
$a_i$ with $a_{i+1}$ and for
$$
i,j\in\{1,\ldots ,n-1\}\ ,\ i<j\quad ,\quad
\sigma_i\cap\sigma_j =
    \cases \emptyset \ \ &\text{if } |i-j|\ge 2\\
           a_{i+1} \ \ &\text{if } j=i+1\,.
    \endcases
$$
Let $H_i = H(\sigma_i)$\,.
We call the ordered system of   half-twists $(H_1,\ldots ,H_{n-1})$ a
frame of $B_n[D,K]$ defined by $(\sigma_1,\ldots ,\sigma_{n-1})$\,, or a frame
of $B_n[D,K]$ for short. \enddefinition

\definition{2.4 Notation}

$[A, B] = ABA^{-1} B^{-1}$.

$\la A, B\ra = ABAB^{-1}A^{-1}B^{-1}$.

$(A)_B= B^{-1} AB$.\enddefinition

\newpage

\proclaim{2.5 Theorem}\ \rom{(E.~Artin's braid group presentation)}
\quad Let $\{H_i\}$ be a frame of $B_n.$ Then
$B_n$
is generated by the half-twists $H_i$   and all the relations
between $H_1 ,\dots, H_{n-1}$ follow from $$
\align [H_i,H_j] &=1\qquad \text {if} \quad |i-j|>1,\\
\langle H_i, H_j \rangle &=1 \qquad \text {if} \quad |i-j|=1,
\\ &1 \leq i,j \leq n-1.
\endalign
$$\ep
\pf \ \cite{A} (or \cite{MoTe4}, Chapter 5).\edm

\proclaim{2.6 Theorem}  Let $\{H_i\}$ be a frame of $B_n.$
Then\roster\item"(i)" for $n\ge 2$,\
$\Center B_n$ is isomorphic to
$ Z$ with a generator
$\Dl_n^2=$\newline $(H_1\cdot\dots\cdot H_{n-1})^n.$
\item"(ii)" $B_2\simeq \BZ$ with a generator $H_1.$\endroster
\endproclaim
\pf \ \cite{MoTe4}, Corollary V.2.3.\edm

\proclaim{2.7 Proposition} There is a natural defined homomorphism $B_n\to S_n$
(symmetric group on $n$ elements) defined by $H_i\to (i\ i+1).$\ep

\pf Since the transpositions $\alpha_i=(i\ i+1)$ satisfy the relations from
Artin's theorem (2.5), the above homomorphism is well defined.\edm

\definition{2.8 Definition} \  $\underline{P_n}$.

 The kernel of the above homomorphism is denoted by
$P_n.$\enddefinition

\definition{2.9  Remark} The transpositions $\a_i$ satisfy a relation that
$H_i$ do not satisfy, which is $\a_i^2=1.$
In fact it is true for any transposition. Under the above homomorphism the
image
of any half-twist is a transposition and  thus any square of a half-twist
belongs to $\ker(B_n\to S_n)$ which is
$P_n.$\enddefinition

\definition{2.10 Definition} \  \underbar{Transversal half-twists, adjacent
half-twist, disjoint half-twist}.

Let $\sigma_1$ and $\sigma_2$ be 2 paths in $D$ with endpoints in $K$ which
do not
intersect $K$ otherwise (like in 2.2). The half-twists $H(\s_1)$ and
$H(\s_2)$ will be
called {\it transversal} if
$\s_1$ and $\s_2$ intersect transversally in one point which is not an end
point of either of the $\s_i$'s.

The half-twists $H(\s_1)$ and $H(\s_2)$ will be called {\it adjacent} if $\s_1$
and $\s_2$ have one endpoint in common.

The half-twists $H(\s_1)$ and $H(\s_2)$ will be called {\it disjoint} if $\s_1$
and $\s_2$ do not intersect.\enddefinition

\proclaim{2.12 Claim} Disjoint half-twists commute and adjacent half-twists
satisfy the
triple relation $ABA=BAB.$\ep

\demo{Proof} By Proposition 2.7 and the fact that every 2 half-twists are
conjugated to
each other.\edm

\definition{2.12 Definition} \  $\underline{\tB_n}.$

Let $Q_n$ be the subgroup of $B_n$ normally generated by $[X,Y]$ for $X,Y$
transversal half-twists.  $\tB_n$ is the quotient of $B_n$ modulo $Q_n.$
For $X\in B_n$ we denote by $\tX$ the image of $X$ in $\tB_n.$\ $\{\tilde
H_i\}$
is a frame of $\tB_n$ if $\{H_i\}$ is a frame of $B_n.$\enddefinition

Later we shall need some basic relations satisfied in $\tilde B_n$ (and not in
$B_n$).
We formulate this in the following claim.

\proclaim{2.13 Claim} Let $\tilde P_n$ be the image of $P_n$
(from $2.7$) in $\tilde B_n.$
Then $\tilde P_n'=\{1,c\}$ where $c^2=1,$\ $c\in\Center\tilde B_n.$
In particular, if $\tX$ and $\tY$ are $2 $ adjacent half-twists
$[\tX^{\pm2},\tY^{\pm2}]=c.$\ep
\pf \ \cite{MoTe9}, Proposition II.5.2.\enddemo

\bk

\subheading{\S3. General results on fundamental groups of complements of
curves}

In this section we prove two general results based on the situation described
in 1.4.
The first one concerns the action of the braid group on the fundamental group
$\pi_1(\BC^2-S)$  and the second one is a corollary on the structure of
$\pi_1(\BC^2-S).$

\proclaim{3.1 Proposition}
Consider the situation of $1.4.$
Let $\Dl_m^2$ be the generator of the center of $B_m[\BC_u,\BC_u\cap S].$
Then  when considered as elements of $\pi_1(\BC^2-S),$\newline
$\Dl_m^2(\G_k)=\G_k$\
\ $\forall \G_k.$\ep

\pf Let $\vp_u$ be the naturally defined homomorphism from
$\pi_1(\BC-N,u)\to$\linebreak $ B_m[\BC_u,\BC_u\cap S].$
This homomorphism is called the braid monodromy and it factors through the
classical monodromy from $\pi_1$ to $ S_m$,\
$
\underbrace{\pi_1\overset\vp_u\to\rightarrow
B_m\to
S_m}_{\psi}.$
Since $B_m$ acts on $\pi_1(\BC_u-S,*),$ so does $\vp_u(\pi_1) $.
Moreover, $\vp_u(\pi_1)$ acts on the elements of $\pi_1(C_u-S,*)$ when
considered as elements of $\pi_1(\BC^2-S,*).$
By the affine Van Kampen theorem (see \cite{RoTe}), for $\g\in\pi_1(\BC-N,u),$\
$\vp_u(\g)(\G_k)=\G_k$ when considered as
elements of $\pi_1(\BC^2-S).$
By \cite{MoTe4}, Lemma VI.2.1, \ $\Dl_m^2$ is a product of elements of the
form $\vp_u(\g )$ for $\g\in\pi_1(\BC-N,u)$ and thus $\Dl_m^2(\G_k)=\G_k\
\forall \G_k.$\enddemo

\proclaim{3.2 Proposition}
Consider the situation of $1.4.$
Let $\G=\prod\limits_{i=1}^m\G_i.$
Then when $\G$ is considered as an element of $\pi_1(\BC^2-S,*),$ it is a
central element and $\langle \G\rangle$ is an infinite cyclic group.\ep

\pf We first consider   $\G=\prod\limits_{i=1}^m
\G_i$ as an element of $\pi_1(\BC_u-S,*).$
Clearly, $\prod \G_i$ is homotopic to a loop $\partial D$ around all the points
of $\BC_u\cap S.$
By Proposition V.2.1 of \cite{MoTe5} in
$\pi_1(\BC_u-S,*)$ the conjugation of $\G_k$ by $\partial D$  is equal to the
action of
$\Dl_m^2$ (which is defined in 2.6) on $\G_k.$
But in $\pi_1(\BC^2-S,*)$,\ $\Dl_m^2$ acts trivially on $\G_k,$ (by 3.1)
so conjugation of $\G_k$ by $\partial D$ in $\pi_1(\BC^2-S)$ is stable and thus
$\G=\partial D$ is in the center of $\pi_1(\BC^2-S,u).$
Moreover, since $\langle\Dl_{m}^2\rangle$ is an infinite cyclic group (see
(\cite{MoTe4}, V.2.1)), so is
$\langle \G\rangle.$\enddemo

\bk

\subheading{\S4. Results on $\pi_1(\BC^2-S,*)$ for a Veronese branch curve}

In this section we restrict ourselves to a curve which is the branch curve of a
Veronese generic projection.
We will quote results concerning its complement in $\BC^2$ (cf. \cite{MoTe9}
and \cite{MoTe10}) which will be used later in the proof of the main result,
concerning  its complement in
$\CPt.$

The fundamental group of the complement in $\BC^2$ turned out to be a
quotient of
a semidirect product of
$\tilde B_{n^2}$ (for a Veronese embedding of $\deg n$) and $G_0(n^2)$ which is
a $\BZ_2$ extension  of a free group on $n^2-1$ elements (see
\cite{MoTe9}).

From now on we will restrict ourselves to a Veronese embedding of $\deg 3.$
Let $S$ be the branch curve of a generic projection to $\CPt$ of a Veronese
surface of $\deg 3.$
The degree of the projection is $9,$
and the degree of the branch curve is 18.
Let $\BC^2$ be a big affine piece of $\CPt$ s.t. $S$ is transversal to the line
in infinity.

Consider $\tB_9$ as defined in \S2.
Instead of working with a frame of $\tB_9$ we will work with $\{\tT_i\},$ a set
of generators for
$\tB_9$ as follows:

\definition{4.1 Definition} Let $\{T_i\}_{i=1\ i\ne 4}^9$  be s.t. $\tT_i$ is a
half-twist w.r.t. $t_i$ where $t_i$ are arranged as in Figure 4.1

\midinsert
\centerline{
\epsfysize=1.5in
\epsfbox{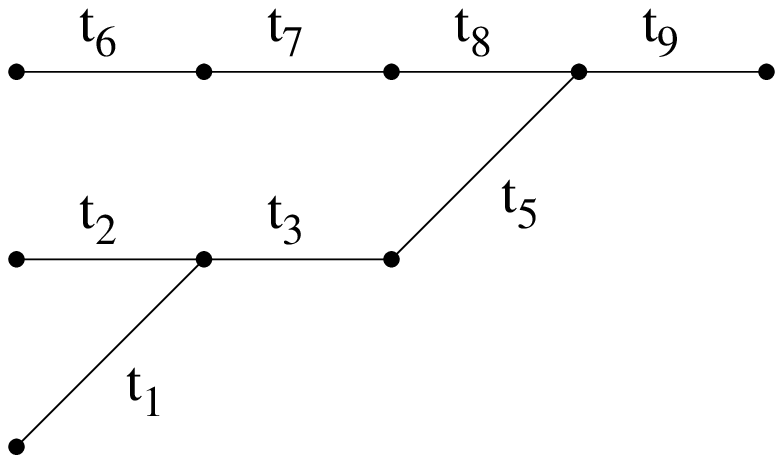}}  
 \botcaption{Figure 4.1}
 \endcaption
\endinsert

\enddefinition

\definition{4.2 Remark}\ The choice of the base originates from a
configuration of planes in the degeneration of $V_3$ to a union of planes.
 We constructed this degeneration in  \cite{MoTe7}, but we do not
use it directly in this paper.  It was used in \cite{MoTe9} to
prove the results which are quoted here. \enddefinition

\subheading{4.3}
 \ It is easy to see that

$T_i$ and $T_j$ are adjacent for $(i,j)$ as follows:

$\qquad \quad i,j \in \{1, 2, 3\}$

$\qquad \quad i=5\qquad \qquad j=3,8,9$

$\qquad\quad i=6,7,8 \qquad j=i+1.$

$T_i$ and $T_j$ are disjoint for $(i,j)$ as follows:

$\qquad\quad i\in \{1,2,3\}\ \ \ j\in \{6,7,8,9\}$

$\qquad\quad i=5\qquad \quad\ \ j= 1,2,6,7$

$\qquad\quad i=6\qquad \quad\ \ j= 8,9$

$\qquad\quad i=7\qquad \quad\ \  j=9.$

\proclaim{4.4 Claim} The set $\{T_i\}$ satisfies the following relations:
$$\alignat 2
& \langle T_i,T_j\rangle=1\quad &&\text{if}\ T_i\ \text{and}\ T_j  \ \text{are
adjacent}\\
&[T_i,T_j]=1\quad &&\text{if}\ T_i\ \text{and}\ T_j  \ \text{are
disjoint}\\
&[T_1,T_2\1T_3T_2]=1\\
&[T_5,T_8\1T_9T_8]=1\endalignat$$
\ep

\pf Since the sequence of  half-twist $\left\{T_1, T_2, T_2\1 T_3 T_2, T_5,
T_9, T_8\1\
T_9\ T_8, T_7, T_6\right\}$ is represented by a consecutive sequence of
paths (see Fig.
4.2),   it is a frame.  By E. Artin's Theorem, they satisfy the relations
that a
frame satisfies (Theorem 2.5). When writing down the triple relations for
the above frame,
we get
 $$\align
&\la T_1,T_2\ra =1\\
&\la T_2,T_2\1T_3T_2\ra=1\\
&\la T_2\1T_3T_2,T_5\ra=1\\
&\la T_5,T_9\ra=1\\
&\la T_9,T_8\1T_9T_8\ra=1\\
&\la T_8\1T_9T_8,T_7\ra=1\\
&\la T_7,T_6\ra=1\endalign$$
When writing down the commutative relations, we get:
$$\align
&[T_i,T_j]=1\ \text{for $T_i$ or $T_j$ disjoint plus}\\
&[T_1,T_2\1T_3T_2]=1\\
&[T_5,T_8\1T_9T_8]=1.\endalign$$
We just need to show $\la T_8,T_9\ra=1,$\ $\la T_2,T_3\ra=1.$
Since $T_8$ and $T_9$ are adjacent, by Claim 2.11, $T_8\1T_9T_8=T_9T_8T_9\1.$
Now from $\la T_9, T_8\1T_9T_8\ra=1,$ we get
$$\align
1&=T_9T_8\1T_9T_8T_9T_9T_8\1T_9\1T_9T_9\1T_8\1T_9\1=\\
&  =T_9T_9T_8T_9\1T_9T_9T_8\1T_9\1T_8\1T_9\1\\
&  = T_9T_9T_8T_9T_8\1T_9\1T_8\1T_9\1
\endalign$$
Thus also $T_9T_8T_9T_8\1T_9\1T_8\1=1.$
Thus $\la T_9,T_8\ra=1.$
Similarly, $\la T_2,T_3\ra=1$, and we get the claim.
\enddemo

\midinsert
\centerline{
\epsfysize=1.5in
\epsfbox{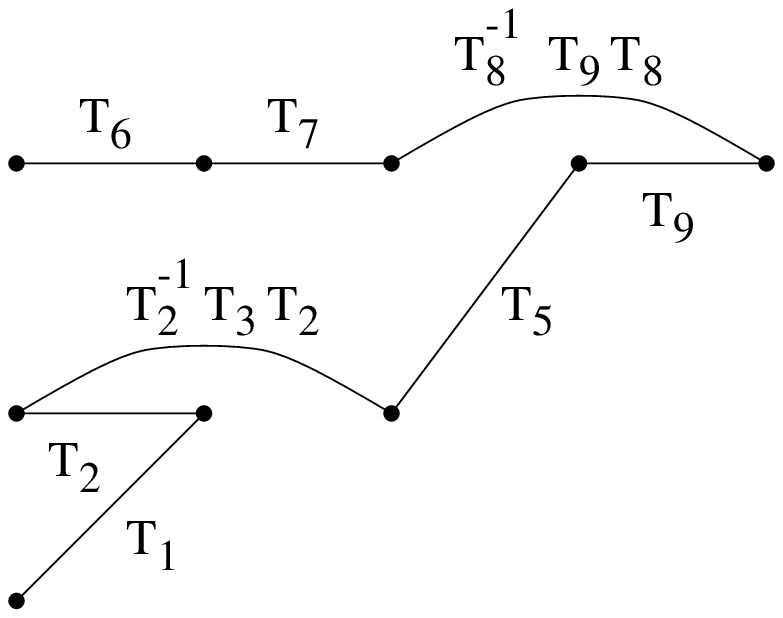}}  
 \botcaption{Figure 4.2}
 \endcaption
\endinsert

\definition{4.5 Definition} \   \underbar{Polarization, orderly adjacent,
non orderly adjacent}.

 We choose an orientation on each $T_i$ with compatibility
with its ``bigger'' neighbor.
We call it a polarization.  See Figure 4.3.

Most of the adjacent $T_i$'s are orderly adjacent (compatible polarization)
apart from $\{T_1,T_2\}$ and $\{T_5,T_8\}$ which are non orderly adjacent.

\midinsert
\centerline{
\epsfysize=1.5in
\epsfbox{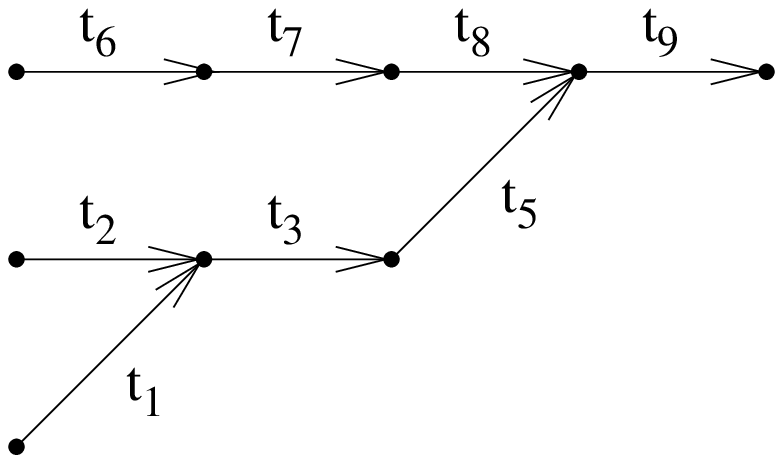}}  
 \botcaption{Figure 4.3}
 \endcaption
\endinsert

\enddefinition

\definition{4.6 Definition} \  $\underline{G_0(9)}.$

$G_0(9)$ is a $\BZ_2$ extension of a free group on 8 elements.
We take the following model for $G_0(9)\:$

Let $G_0(9)$ be generated by $\{g_i\}_{i=1\ i\ne4}^9$ s.t. $$[g_i,g_j]=\cases
1\quad&T_i, T_j\ \text{are disjoint}\\
\tau\quad & \text{otherwise}\endcases$$
where $\tau^2=1,$\ $\tau\in\Center G_0(9).$

We take the following action of $\tB_9$ on $G_0(9)$
$$(g_i)_{\tT_k}=\cases g_i\1\tau&\quad k=i\\
g_i&\quad T_i,T_k\ \text{are disjoint}\\
g_ig_k\1&\quad T_i,T_k\ \text{are not orderly adjacent}\\
g_kg_i&\quad \text{otherwise}\endcases$$\enddefinition

\definition{4.7 Definition}\  $\underline{G_ 9,\ c}.$

Consider the semidirect product $\tB_9\ltimes G_0(9)$ w.r.t. the chosen action.

Let $c=[\tT_1^2,\tT_2^2].$

Let $\xi_1=(\tT_2\tT_1\tT_2\1)^2\tT_2^{-2}.$

Let $N_9\triangleleft\tB_9\ltimes G_0(9)$ be normally generated by $c\tau\1$
and $(g_1\xi_1\1)^3.$

Let $G_9=\df{\tB_9\ltimes G_0(9)}{N_9}.$\enddefinition

\definition{4.8 Definition}\  $\underline{\hat\psi_9}.$

Let $\tilde\psi_9$ be the homomorphism $\tB_9\to S_9$ induced from the standard
homomorphism $B_9\to S_9$ (see 2.7).
$\tilde\psi_9$ exists since $[X,Y]\to1$ under the standard homomorphism.
Let $\hat\psi_9: G_9\to S_9$ be defined by the first coordinate
$\hat\psi_9(\a,\be)=\tilde\psi_9(\a).$\enddefinition

\definition{4.9 Definition}\  $\underline{\psi }.$

The projection $V_3\to \BC^2$, of degree 9, induces a standard monodromy
homomorphism
 $\pi_1(\BC^2-S,*)\to S_9$ which we denote by $\psi .$\enddefinition

\proclaim{4.10 Proposition}
$\pi_1(\BC^2-S,*)\simeq G_9$ s.t. $\psi$ is compatible with $\hat\psi_9$.\ep

\pf \ \cite{MoTe9}, VI.1.\enddemo

\definition{4.11 Definition}\  $\underline{H_9,H_{9,0},H_9',H_{9,0}'}.$

Let $Ab: B_9\to\Bbb Z$ be the abelianization of $B_9$  and $B_9$ over its
commutator
subgroup.

Let $\widetilde{A_b}:\tB_9\to\BZ$ be a homomorphism induced from $Ab$
(which exists
since $Ab([X,Y])=1).$

Let $\widehat{Ab}: G_9\to\BZ$ be defined by the first coordinate
$\widehat{Ab}(\a,\be)=\widetilde{Ab}(\a).$

Let $H_9=\ker\hat\psi_9.$

Let $H_{9,0}=\ker \hat\psi_9\cap\ker \widehat{Ab}.$

Let $H_9',H_{9,0}'$ be the commutant subgroup of $H_9$ and $H_{9,0}$
respectively.\enddefinition

\proclaim{4.12 Proposition}
There exists a series
$1\triangleleft H_{9,0}'\triangleleft H_{9,0}\triangleleft H_9\triangleleft
G_9,$
where $G_9/H_9\simeq S_9,$ \quad\qquad $H_9/H_{9,0}\simeq\Bbb Z,$\quad
\qquad  $H_{9,0}/H'_{9,0}\simeq
(\Bbb Z\oplus\Bbb Z/3\BZ)^8,$\qquad $H_{9,0}'=H_9' \cong \Bbb Z/2 \Bbb Z.$
\endproclaim

\demo{Proof} \ \cite{MoTe10}, Proposition 2.4.\enddemo

Our main result is a result of type 4.12.

To this end we need to get into the proofs of the structure theorems for
$G_9,$ which are quoted in  4.10
and 4.12.  We need this for the proof of the main result

\subheading{4.13}\
$H_{9,0}$ is generated by
 $\{g_i\}_{i=1\ i\ne 4}^9, \quad \{\xi_i\}_{i=1\ i\ne 4}^9,\quad c$
 where
$$[g_i,g_j]=\cases 1&\quad T_i,T_j\ \text{are disjoint}\\
c&\quad \text{otherwise}.\endcases$$
$$[\xi_i,\xi_j]=\cases         1&\quad T_i,T_j\ \text{are disjoint}\\
c&\quad \text{otherwise}.\endcases$$

$$[\xi_i,g_j]=\cases         1&\quad T_i,T_j\ \text{are disjoint}\\
c&\quad \text{otherwise}.\endcases$$

$g_i,\xi_i$ of infinite order.

$g_i^3=\xi_i^3.$

$c^2=1.$

$c\in\Center (G_9).$

$H_9$ is generated by $H_{9,0}$ and $\tilde T_1^2$ where $\tilde T_1^2$ is of
infinite order.

$H_9'=H_{9,0}'$ is generated by $c.$
($c$ is the image of the generator of $\tilde P'_9$ from 2.12).

For simplicity we also denote $\zeta _i=g_i\xi_i\1,$ ($\zeta_i^3=1).$

\medskip

From 4.13 we get the following:

\subheading{4.14}\   $H_{9,0}'=H_9'\simeq\BZ_2 (\subseteq\Center(G_9)).$

$H_{9,0}/H_{9,0}'$
is generated by $\{\xi_i\}_{i=1,i\ne 4}^9$ and $\{\zeta _i\}_{i=1,i\ne 4}^9$,
when the only
relations are the commutativity relation and $\zeta _i^3=1.$
Thus $\df{H_{9,0}}{H_{9,0}'}\simeq (\BZ\oplus\BZ/3\BZ)^8.$

$H_9/H_{9,0}$ is generated by $\tT_1^2$ and thus is isomorphic to $\BZ.$

  $H_9$ is the kernel of $G_9\to S_9,$ and thus $G_9/H_9\simeq S_9.$

\bk

\subheading{\S5.\ The Main Result}

Our main result is the following theorem.

\proclaim{5.0 Theorem}
Let $S$ be the branch curve of a generic projection to $\CPt$ of a Veronese
embedding of $\deg 3.$
Then $\pi_1(\CPt-S,*)$ is an extension of a solvable group by the symmetric
group of $9$ elements.
In fact, we have $1\triangleleft \ov H_{9,0}'\triangleleft \ov
H_{9,0}\triangleleft \ov
H_9\triangleleft  \ov G_9$ where $\ov G_9/\ov H_9\simeq S_9,$\ $\ov H_9/\ov
H_{9,0}\simeq
Z_9,$\ $\ov H_{9,0}/ \ov H_{9,0}' \simeq(\BZ\oplus\BZ_3)^8$, $\ov
H_{9,0}'\simeq \BZ_2.$\ep

\pf
It is easy to calculate $\deg S$ (see \cite{MoTe3}) and it is 18.

We consider the situation of 1.4 for the branch curve from our Theorem.
By \cite{MoTe9}, Lemma 2.3, there is a possibility to choose a
$g$-base
$\{\G_i,\G_{i'}\}_{i=1}^9 $  s.t. $\psi(\G_i)=\psi(\G_{i'})=$ transposition.
(This choice is    a consequence of the
degeneration of the surface to a union of 9 planes.)

By 1.5, \
$\pi_1(\CPt-S,*)=\df{\pi_1(\BC^2-S,*)}{\big
\langle\prod\limits_{i=9}^1\G_{i'}\G_{i}\big\rangle}.$
 Denote
$\hat\be: \pi_1(\CPt-S,*)\to G_9$ to be the isomorphism
from 4.10 and
$\dl=\hat\be\left(\prod\limits_{i=9}^1\G_{i'}\G_{i}\right).$
Clearly, $\pi_1(\CPt-S,*)\simeq \df{G_9}{\langle\dl\rangle}$ which we denote by
$\ov G_9.$
To prove the theorem we shall prove that $\ov G_9=\df{G_9}{\langle \dl\rangle}$
is an extension of a solvable group by a symmetric group.

In 4.12 we introduced  a sequence $1\triangleleft H_{9,0}'\triangleleft
H_{9,0}\triangleleft H_9\triangleleft G_9$ and the appropriate quotients.
Let $\ov H_{9,0}',\ \ov H_{9,0},\ \ov H_9$ be the images of $H_{9,0}',
\ H_{9,0}, H_9$ in $\ov G_9$ respectively.
To prove the theorem we shall compute  $\ov G_9/\ov H_9,$\ $\ov H_9/\ov
H_{9,0},$\
$\ov H_{9,0}/\ov  H_{9,0}$ and $\ov H_{9,0}'.$

We first need to prove some results on $G_9$ in general and on $\dl$ in
particular.
This is done in Claims 5.1--5.10.
From general arguments we already know that $\dl\in\Center(G_9)$ (cf.
Proposition 3.2).\enddemo

\proclaim{5.1 Claim} \ $\dl\in H_9,$\ $\hat\be(\G_i\G_{i'})\in H_9.$\ep

\demo{Proof of Claim 5.1}
By 4.10, $\hat\psi_9\hat\be(\G_{i'}\G_i)=\psi_ (\G_{i'}\G_i).$
Since $\psi (\G_i)=\psi (\G_{i'})=$ transposition.
 $\psi (\G_{i'}\G_i)=1,$ and  thus           $\hat
\be(\G_{i'}\G_{i})\in\ker\hat\psi_9 =H_9.$
Since $\dl=\prod\hat\be(\G_{i'}\G_{i}),$ it is also in
$H_9.$\hfill$\qed$ for Claim 5.1\enddemo

The new quotients will be determined by an expression of $\dl$ as a product of
elements in $H_{9,0}$ and elements which are in $H_9$ but not in $H_{9,0}.$

\definition{5.2 Definition}

By abuse of notation the images in  $G_9=\df{\tB_9\ltimes  G_0(9)}{N_9},$ of
$\tT_i$ from $\tB_9$ (see 4.2)  are also denoted by $\tT_i.$
We also define:

$\tT_4=(\tT_5)_{\tT_8\1\tT_7\tT_3\1\tT_2} $.

$g_4=(g_5)_{\tT_8\1\tT_7\tT_3\1\tT_2} $ for $g_5$ from 4.13.

$\xi_4=(\xi_5)_{\tT_8\1\tT_7\tT_3\1\tT_2} $ for $\xi_5$   from
4.13.\enddefinition

To work in $G_9$ we need some commutativity relations:
\proclaim{5.3 Claim} In $G_9$: \roster\item"(i)" $\tT_i^2,$\
  $\xi_{i},g_i \in H_9 $\  $i=1,\dots,9.$
\item"(ii)" $[\tT_i^2,g_j],[\xi_ig_i]= 1$ or $c.$
\item"(iii)" If $X$ and $Y$ are 2 adjacent half-twists, then
$[\tX^2,\tY^2]=c.$\endroster\ep
\demo{Proof}

(i)\
 Since $ T_i$ is a half-twist $i=1\dots 9,$ thus $\hat\psi_9(\tT_i)$ is a
transposition and $\hat\psi_9(\tT_i^2)=1.$
Thus $\tT_i^2\in H_9.$
The elements $\{g_i,\xi_i\}_{i=1\ i\ne 4}^9$ are in $H_9$ by 4.13.
Since $g_5$ is in $H_9$  and $H_9$ is a normal subgroup $(=\ker \hat
\psi_9),$ \
$g_4$ is also in $H_9.$ The same applies for $\xi_4.$

(ii)\ Since $H_9'=\{1,c\},$\ $c^2=1.$

(iii)\ Since it is true in $\tB_9$ by Claim 2.13.
\hfill$\qed$ for Claim 5.3\edm

In order to compute the corresponding quotients in $\ov G_9$, we need to
express
$\dl$ which is in $H_9  $ (see 5.1) in terms of the following generators of
$H_9$:
$\{\zeta_i\}_{i=1\ i\ne 4}^9,$\ $\{\xi_i\}_{i=1\ i\ne 4}^9,c$ and $\tT_1^2 $
(see 4.13).
Recall that   $\dl=\prod\limits_{i=9}^1\hat\be(\G_{i'}\G_i)$ where
$\fa i=1\dots 9$
$\hat\be(\G_{i'}\G_i)\in H_9$.
Thus, we shall first express $\hat\be(\G_{i'}\G_i)$ for $i=1\dots 9$ in
terms of
$\{\zeta_i\}_{i=1\ i\ne 4}^9,$\ $\{\xi_i\}_{i=1\ i\ne 4}^9,c$ and
$\tT_1^2$, and
then we multiply these expressions to get an expression for $\dl $ in these
generators (see 5.9).
In 5.10 we replace $g_i$ by $\zeta_i\xi_i.$
\proclaim{5.4 Claim}
$$\align
& \hat\be(\G_{1'}\G_1)=g_1\tT_1^2\quad \\
&\hat\be(\G_{2'}\G_2)=g_2\1\xi_2\tT_2^2 \\
&\hat\be(\G_{3'}\G_3)=g_3\xi_3\1\tT_3^2\\
&\hat\be(\G_{4'}\G_4)=g_4\1\xi_4\tT_4^2 \\
&\hat\be(\G_{5'}\G_5)=g_5\1\xi_5\tT_5^2\\
&\hat\be(\G_{6'}\G_6)=g_6\tT_6^2\\
&
\hat\be(\G_{7'}\G_7)=g_7\xi_7\1\tT_7^2\\
&\hat\be(\G_{8'}\G_8)=g_8\1\xi_8\tT_8^2\\
&\hat\be(\G_{9'}\G_9)=cg_g\1\tT_9^2\endalign $$\ep

\pf We take a new  set of generators for $G$:
$$E_i=\cases \G_i\quad & i\ne 2,7\\
\G_{i'}\quad &i=2,7\endcases\quad
E_{i'}'=\cases \G_{i'}\quad & i\ne 2,7\\
\G_{i'}\G_i\G_{i'}\1\quad& i=2,7\endcases$$

(This choice which was made in \cite{MoTe9} originated from a certain
relation in
$G$ induced by the affine Van Kampen Theorem.)
Clearly $\G_{i'}\G_i=E_{i'}E_i.$
Let $A_i=E_{i'}E_i\1.$
Clearly, $\G_{i'}\G_i=E_{i'}E_i=A_iE_i^2.$
By the construction of $\hat\be $ (see \cite{MoTe9}, Ch.V), $
\hat\be(E_i^2)=\tT_i^2$ and $\hat\be (A_i)$ is as follows:
$$\align&\be(A_1)=g_1\\
&\be(A_2)=g_2\1\xi_2\\
&\be(A_3)=g_3\xi_3\1\\
&\be(A_4)=g_4\1\xi_4\\
&\be(A_5)=g_5\1\xi_5\\
&\be(A_6)=g_6\\
&\be(A_7)=g_7 \xi_7\1\\
&\be(A_8)=g_8\1\xi_8\\
&\be(A_9)=cg_9\1 \endalign$$
\hfill $\qed$ for Claim 5.6\edm

In the next step we express $\tT_i^2$ in terms of $\{\xi_i\}_{i=1}^9$ and
$\tT_1^2.$
The main point in the proof of the next claim is that for 3 half-twists
which form
a triangle where one of the edges is $T_i$, the product $\tX^2\tY^{-2}$ of the
other 2 half-twists can be expressed in terms of $\xi_i,\xi_i\1$ and $c.$
The exact statement is as follows:

\proclaim{5.5 Claim} Let $X,Y$ be $2$ half-twists,\quad $X=H(x)$,\quad
$Y=H(y),$\quad $T_i=H(t_i)$ s.t. $x,y,t_i$ make a triangle.
Assume that $x$ and $y$ meet in $\nu,$ and a counterclockwise rotation around
$\nu$ inside the triangle meets $x$ before it meets $y.$
\roster\item"(i)" If the polarization of $T_i$ goes from $x$ to $y,$ then
$\xi_i=\tX^2\tY^{-2}.$
\item"(ii)" If the polarization of $T_i$ goes from $y$ to $x,$ then
$\xi_i=\tX^{-2}\tY^2.$\endroster\ep
\pf Claim IV.4.1 of \cite{MoTe9}.\edm

\newpage

Using this claim we prove the following
\proclaim{5.6 Lemma}
\roster\item"(i)" $\tT_2^2\tT_1^{-2}=c\xi_1\1\xi_2$
\item"(ii)" $\tT_3^2\tT_1^{-2}=\xi_1\1\xi_3\1$
\item"(iii)" $\tT_5^2\tT_1^{-2}=\xi_1\1\xi_3^{-2}\xi_5\1$
\item"(iv)" $\tT_9^2\tT_1^{-2}=c\xi_1\1\xi_3^{-2}\xi_5^{-2}\xi_9\1$
\item"(v)" $\tT_8^2\tT_1^{-2}=\xi_1\1\xi_3^{-2}\xi_5^{-2}\xi_8$
\item"(vi)"
$\tT_7^2\tT_1^{-2}=c\xi_1\1\xi_3^{-2}\xi_5^{-2}\xi_7\xi_8^2$
\item"(vii)"
$\tT_6^2\tT_1^{-2}=\xi_1\1\xi_3^{-2}\xi_5^{-2}\xi_6\xi_7^2\xi_8^2$
\item"(viii)"
$\tT_4^2\tT_1^{-2}=\xi_1\1\xi_3^{-2}\xi_4\xi_5^{-2}\xi_7^2\xi_8^2$\endroster\ep

\pf The proof is based on Claim 5.5.

Moreover, we interchange between the $\xi_i$'s using the commutator from 4.13
$$[\xi_i^{\pm1},\xi_j^{\pm1}]=\cases 1&\quad T_i\ \text{and}\ T_j\ \text{are
disjoint}\\ c&\quad \text{otherwise}\endcases$$

(i) \ We write $T_2^2\tT_1^{-2} =\tT_2^2 \tilde{W}^{-2}\tilde{W}^2\tT_1^{-2}$
for $W=(T_1)_{T_2\1}$ which creates a triangle with $\tT_1$ and $\tT_2$ (see
Fig. 5.1)
We use Claim 5.5 twice --
first when we take $W,T_2,T_1$ instead of $X,Y,T $ from Claim 5.5(i), and
second when we take $T_1,W,T_2$ instead of $X,Y,T_2$ from Claim 5.5(ii). By
Claim 5.5(i)
$\tilde W^2\tT_2^{-2}=\xi_1,$ and thus
$\tT_2^2\tilde W^{-2}=\xi_1\1.
$   By Claim 5.5(ii) $\tT_1^{-2}\tilde W^2=\xi_2,$ and since by Claim 5.3(iii)
$[\tT_1^{-2},\tilde{W}^2]=c,$ we get
$\tilde{W}^2\tT_1^{-2}=c\xi_2.$
Together we get (i).

\midinsert
\centerline{
\epsfysize=1in
\epsfbox{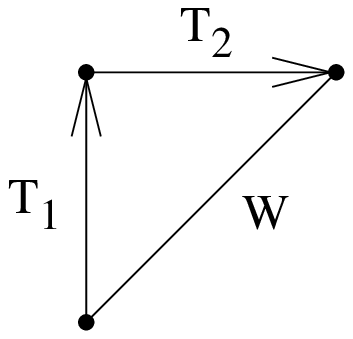}}  
 \botcaption{Figure 5.1}
 \endcaption
\endinsert

(ii)   We write
$\tT_3^2\tT_1^{-2}=(\tT_3^2\tZ^{-2})(\tZ^2\tT_1^{-2})$ for
$Z=(\tT_1)_{\tT_3},$ which creates a triangle with $T_1$ and $T_3.$ (See
Fig. 5.2)
By Claim 5.5 applied twice, $\tT_3^{-2}\tZ^2=\xi_1$ and
$\tZ^{-2}\tT_1^2=\xi_3.$
Thus
$\tZ^{-2}\tT_3^2= \xi_1\1$ and
$\tT_1^{-2}\tZ^2= \xi_3\1.$ By 5.3(iii) we get
$\tT_3^2\tZ^{-2}=c\xi_1\1,$\ $\tZ^2\tT_1^{-2}=c\xi_3\1$.
Since
$c\in\Center(G_9)$ and $c^2=1,$ then $\tT_3^2\tT_1^{-2}=\xi_1\1\xi_3\1.$

\midinsert
\centerline{
\epsfysize=1in
\epsfbox{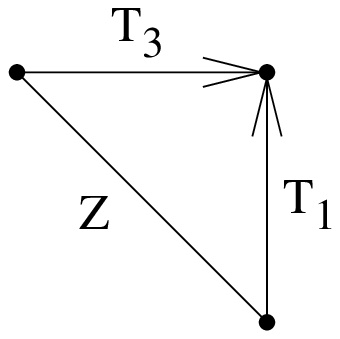}}  
 \botcaption{Figure 5.2}
 \endcaption
\endinsert

(iii) We install $\tT_3^{-2}\tT_3^{2}$ in the middle.
Since $\tT_3$ relates to
$\tT_5$ as
$\tT_1$ relates to
$\tT_3$ we have
$\tT_5^2\tT_3^{-2}=\xi_3\1\xi_5\1.$
Then $\tT_5^2\tT_3^{-2}\tT_3^2\tT_1^{-2}=\xi_3\1\xi_5\1\xi_1\1\xi_3\1.$ Since
$[\xi_3,\xi_5]=[\xi_3,\xi_1]=c,$ \  $[\xi_1,\xi_5]=1,$  \ $c^2=1,$ and
$c\in\Center(G_9),$ we get
$\tT_5^2\tT_1^{-2}=\xi_1\1\xi_3^{-2}\xi_5\1.$

(iv) We install $\tT_5^{-2}\tT_5^2$ in the middle.
Since $\tT_9$ relates to $\tT_5$
as
$\tT_3$ relates to
$\tT_1,$ then
$\tT_9^2\tT_5^{-2}=\xi_5\1\xi_9\1.$
Thus
$\tT_9^2\tT_1^{-2}=\tT_9^2\tT_5^{-2}\tT_5^2\tT_1^{-2}=\xi_5\1\xi_9\1\xi_1\1\xi_3
^{-2}
\xi_5\1=c\xi_1\1\xi_3^{-2}\xi_5^{-2}\xi_9\1.$
The last equation is based on the commutators of $\xi_i$ and the fact that
$c\in \Center G_9,$\ $c^2=1$.

(v)\ $\tT_8^2\tT_1^{-2}=\tT_8^2\tT_5^{-2}\tT_5^2\tT_1^{-2}.$
\ $T_8$ relates to $T_5$ as $T_2$ relates to $T_1$ and thus
$\tT_8^2\tT_5^{-2}=c\xi_5\1\xi_8.$
Thus $\tT_8^2\tT_1^{-2}=c\xi_5\1\xi_8\xi_1\1\xi_3^{-2}\xi_5\1.$
Since $\xi_8 $ commutes with $\xi_1,\xi_3$ and $\xi_5$ commutes with
$\xi_1$ and
$[\xi_5,\xi_8]=[\xi_5,\xi_3]=c,$ we have $\tT_8^2\tT_1^{-2}  =
c^4\xi_1\1\xi_3^{-2}\xi_5^{-2}\xi_8$ which equals
$\xi_1\1\xi_3^{-2}\xi_5^{-2}\xi_8.$

(vi)\ $\tT_7^2\tT_1^{-2}=\tT_7^2\tT_8^{-2}\tT_8^2\tT_1^{-2}.$
Since $T_7$ relates to $T_8$ as $T_1$ relates to $T_3,$ then
$\tT_8^2T_7^{-2}=\xi_7\1\xi_8\1$ and $\tT_7^2 \tT_8\2 = \left(\xi\1_7
\xi_8\1\right)\1 = \xi_8 \xi_7$ and thus
$\tT_7^2\tT_1^{-2}=\tT_7^2\tT_8^{-2}\tT_8^2\tT_1^{-2}=(\xi_8\xi_7) \cdot
\xi_1\1\xi_3^{-2}\xi_5^{-2}\xi_8.$
As before,
$\tT_7^2\tT_1^{-2}=c^3\xi_1\1\xi_3^{-2}\xi_5^{-2}\xi_7\xi_8^2 $ which
equals\linebreak
$c\xi_1\1\xi_3^{-2}\xi_5^{-2}\xi_7\xi_8^2.$

(vii)\ $\tT_6^2\tT_1^2=\tT_6^2\tT_7^{-2}\tT_7^2\tT_1^{-2}.$\ $T_6$ relates to
$T_7$ as $T_7$ relates to $T_8$ and thus $\tT_6^2\tT_7^{-2}=\xi_7\xi_6$ (see
(vi)).
Therefore
$\tT_6^2\tT_1^2=c\xi_7\xi_6\xi_1\1\xi_3^{-2}\xi_5^{-2}\xi_7\xi_8^2$
which equals as before $c^2\xi_1\1\xi_3^{-2}
\xi_5^{-2}\xi_6\xi_7^2\xi_8^2=\xi_1\1\xi_3^{-2}\xi_5^{-2}\xi_6\xi_7^2\xi_8^2.$

(viii) \ $T_4=(T_5)_{T_8\1T_7T_3\1T_2}.$ By Claim II.1.0 of \cite{MoTe9}, if
$X=H(x)$ is represented by a diffeomorphism $\be$ and $Y=H(y)$, then
$Y_X=H((y)\be).$ Therefore $T_4$ is a half-twist and we write $T_4=H(t_4).$
Moreover, to
describe
$T_4$ we must apply
$T_8\1,T_7,T_3\1,T_2$ on $t_5$ and we get   $t_4$ is as in
 Figure 5.3:
\midinsert
\centerline{
\epsfysize=1.5in
\epsfbox{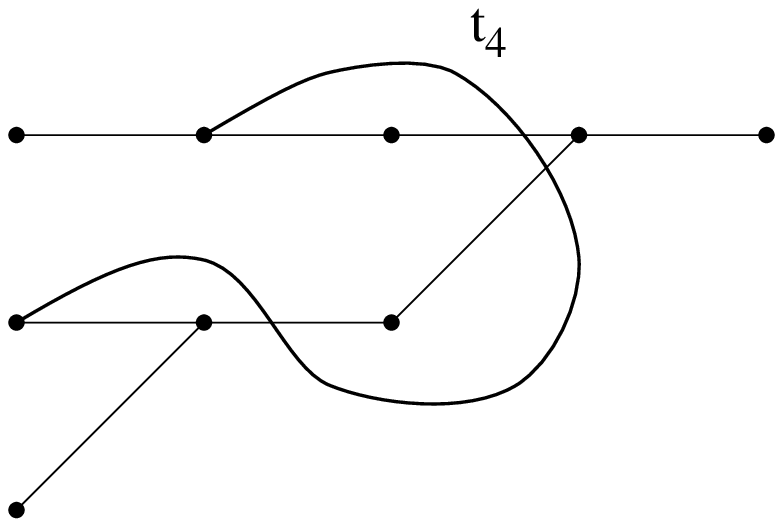}}  
 \botcaption{Figure 5.3}
 \endcaption
\endinsert
Since $T_6$ relates to $T_4$ as $T_2$ relates to $T_1$,then
$\tT_6^2\tT_4^{-2}=c\xi_4\1\xi_6.$
Now
 $\tT_4^2\tT_1^{-2}   =
\tT_4^2\tT_6^{-2}\tT_6^2\tT_1^{-2}=(c\xi_4\1\xi_6)\1\xi_1\1\xi_3^{-2}
\xi_5^{-2}\xi_6\xi_7^2
\xi_8^2$ which equals as before to
$\xi_1\1\xi_3^{-2}\xi_4\xi_5^{-2}\xi_7^2\xi_8^2 .$

\hfill$\qed$ for Lemma 5.6\enddemo

In order to express $g_4$ and $\xi_4$ in terms of $\{g_i\}_{i=1\ i\ne 4}^9$ and
$\{\xi_i\}_{i=1\ i\ne 4}^9,$ we need the following claim from \cite{MoTe9}
\proclaim{5.7 Claim}
For $f_i=g_i$ or $\xi_i$
$$(f_i)_{\tT_k} = \cases
f\1_i \nu \qquad k=i\\
f_i \qquad\ \ \ \ T_i, T_k\qquad \text{weakly disjoint}\\
f_kf_i\qquad\  T_i, T_k\qquad \text{orderly adjacent}\\
f_if_k\1\quad\ \ T_i, T_k\qquad \text{are not orderly adjacent}.\endcases$$
$$(f_i)_{\tT_k\1} = \cases
f\1_i \nu \qquad k=i\\
f_i \qquad\ \ \ \ T_i, T_k\qquad \text{weakly disjoint}\\
f_if_k\qquad\  T_i, T_k\qquad \text{orderly adjacent}\\
f_k\1 f_i\quad\ \ T_i, T_k\qquad \text{are not orderly adjacent}.\endcases$$\ep

\pf \cite{MoTe9}, Lemma IV.6.3.
The conjugations for $g_i$ are part of the definition of $G_0(9)$ (see 4.6) and
remains when moving to $G_9=\df{\tB_9\ltimes G_0(9)}{N_9}.$

\hfill $\qed$ for Claim 5.7\edm

Now we can express $g_4,\xi_4$ in terms of $\{g_i\}^9_{i=1\ i\ne 4}.$

\proclaim{5.8 Claim}
\roster\item"(i)" $\xi_4=c\xi_2\xi_3\xi_5\xi_7\1\xi_8\1$
 \item"(ii)" $g_4=cg_2g_3g_5g_7\1g_8\1$\endroster\ep

\pf The proof is similar for (i) and (ii), and is based on 5.9.
We shall only  prove (i).
By 5.2, $\xi_4=(\xi_5)_{\tT_8\1\tT_7\tT_3\1\tT_2}$.
Since $T_5$ and $T_8$ are not orderly adjacent by 5.7,
$(\xi_5)_{\tT_8\1}=\xi_8\1\xi_5.$
Since $T_5$ and $T_7$ are  s disjoint, $(\xi_5)_{\tT_7}=\xi_5.$
Since $T_7$ and $T_8$ are orderly adjacent, $(\xi_8)_{\tT_7}=\xi_7\xi_8$ and
thus $(\xi_8\1)_{\tT_7}=\xi_8\1\xi_7\1.$
Together we have $(\xi_5)_{\tT_8\1\tT_7}=\xi_8\1\xi_7\1\xi_5.$
We now apply $\tT_3\1.$
Since $T_3$ is disjoint from $T_7$ and $T_8,$ then
$(\xi_8\1\xi_7\1)_{\tT_3\1}=\xi_8\1\xi_7\1.$
On the other hand, $T_5$ and $T_3$ are orderly adjacent and thus
$(\xi_5)_{\tT_3\1}=\xi_5\xi_3 $ and
$(\xi_5)_{\tT_8\1\tT_7\tT_3\1}=\xi_8\1\xi_7\1\xi_5\xi_3.$
Now
  $\tT_2$ acts  on  $\xi_3$ to get
$\xi_2\xi_3$ and does not move  the other factors. Thus
$\xi_4=\xi_8\1\xi_7\1\xi_5\xi_2\xi_3.$ We  rearrange the factors using the
comutators of 4.13 to get
$$\xi_4=c\xi_2\xi_3\xi_5\xi_7\1\xi_8\1.$$
\hfill $\qed$ for Claim 5.8\edm

In fact we are interested in $\dl$ up to a product with $c $
and thus we formulate the following:
\proclaim{5.9 Claim} Up to multiplication by $c$
$$\dl=g_1g_2^{-2}g_5^{-2}g_6g_7^2g_9\1\xi_1^{-8}\xi_2^4\xi_3^{-12}\xi_5^{-8}
\xi_6\xi_7
^2\xi_8^6\xi_9\1\tT_1^{18}.$$  \ep

\demo{Proof of Claim 5.9}
By  definition, $\dl=\hat\be\left(\prod\limits_{i=9}^1\G_{i'}\G_i\right)$
which equals $\prod\limits_{i=9}^1\hat\be(\G_{i'}\G_i).$
We substitute in the product the values of $\hat\be(\G_{i'}\G_i),$\
$i=1,\dots,9$ from Claim 5.4.
In the resulting formula, we replace $\tT_i^2$ by
$\left(\tT_i^2\tT_1^{-2}\right)\tT_1^2$ for each $i=1,\dots,9.$
We then substitute the formula for $\tT_i^2\tT_1^{-2}$ from Claim 5.6.
We also substitute the values of $g_4$ and $\xi_4$  from Claim 5.8.
We then get a formula for $\dl$ as a product of $\{\xi_i,g_i\}_{i=1\ i\ne
4}^9$ and $\tT_1^2.$
Since we are not interested in the appropriate power of $c,$ we can use Claim
5.3(ii) by ``pushing'' all the powers of
$\tT_1^2$ to the right end of the last formula and rearrange the $\xi_i$'s and
the $g_i$'is to the get the claim.

\hfill $\qed$ for Claim 5.9\edm

\proclaim{5.10 Claim} Up to multiplication by $c,$ for $\zeta_i$ and $\xi_i$
from 4.13, we have
$\dl=\zeta_1\zeta_2\zeta_5\zeta_6(\zeta_7\zeta_9)\1\xi_1^{-7}\xi_2^2\xi_3^{-12}\
xi_5
^{-10}\xi_6^2\xi_7^4\xi_8^6\xi_9^{-2}\tT_1^{18}.$\ep
\demo{Proof of Claim 5.10} In the formula from Claim 5.9 we replace $g_i$
by $\zeta_i\xi_i.$
Since $H_9'=\{1,c\}$ where $c\in \Center(G_9),$ we can rearrange the terms of
the formula.
Also using $\zeta_i^3=1$ we get the claim.\hfill$\qed$ for Claim 5.10\edm
\medskip
We go back to the proof the theorem.

To prove that $\ov G_9$ is an extension of a solvable group by a symmetric
group, it is enough to find a normal subgroup whose quotient is $S_9.$ The
subgroup will be $\ov H_9.$

Recall (3.12) that there exists a series $1\triangleleft H_{9,0}'\triangleleft
H_{9,0}\triangleleft H_9\triangleleft G_9.$
We defined $\ov G_9=\df{G_9}{\langle\dl\rangle},$  and $\ov
H_9,$\  $\ov
H_{9,0},$\ $\ov
H_{9,0}'$ to be the images of $H_9,$\ $H_{9,0},$\ $H_{9,0}'$ in $\ov G_9$
respectively, and we have a series
 $1\triangleleft \ov H_{9,0}'\triangleleft
\ov H_{9,0}\triangleleft\ov H_9\triangleleft \ov G_9.$
We shall compute the quotients.
Since $\dl\in H_9$ (Claim 5.1),  $\df{\ov G_9}{\ov H_9} \simeq
\df{G_9}{H_9}\simeq S_9.$
Since $\df{H_9}{H_{9,0}}$ is generated by $\tT_1^2$ (see 4.14), $\df{\ov
H_9}{\ov H_{9,0}}$ is also generated by $\tT_1^2.$
By Claim 5.12, $(\tT_1^2)^9\dl\1\in H_{9,0}.$
So when considered as elements of $\ov H_9,$\ $(\tT_1^2)^9\in\ov H_{9,0},$ and
thus as elements of $\df{\ov H_9}{\ov H_{9,0} },$ \ $\tT_1^2$ is of order 9.
Thus $\df{\ov H_9}{\ov H_{9,0}}\simeq \BZ_9.$
Now let
$Y_1=\xi_1^{-7}\xi_2^2\xi_3^{-12}\xi_5^{-10}\xi_6^2\xi_7^4\xi_8^6\xi_9^{-2}(\tT_
1
)^{18}.$
We complete $Y_1$ to a base $Y_1,\dots ,Y_9$ of $\big \langle
\tT_1^2,\{\xi_i\}_{i=1\ i\ne 4}^9\big\rangle.$\
$\df{H_9}{H_9'}$ is generated by $Y_1,\dots, Y_9,\{\zeta_i\}_{i=1\ i\ne 4}^2$.
Modulo $\langle\dl\rangle,$\ $Y_1\in\{\zeta_i\}_{i=1\ i\ne 4}^9$ and
$Y_2,\dots, Y_8$ are of infinite order.
Thus $\df{\ov H_{9,0}}{\ov H_{9,0}'}\simeq\left(\BZ\oplus \BZ/3\BZ\right)^8.$

By 3.2, $\la\dl\ra$ is an infinite cyclic group where $\dl\in\Center (G_9).$
Since $c^2=1,$ then $\la c\ra\cap\la\dl\ra=1$ and $\ov
H_{9,0}'\simeq\df{H_{9,0}'}{\la\dl\ra \cap H_{9,0}'}\sim H_{9,0}'\simeq \BZ_2.$
Thus we have a series $1\triangleleft \ov H_{9,0}'\triangleleft \ov
H_{9,0}\triangleleft \ov H_9\triangleleft \ov G_9$ s.t.  $\ov H_9$ is a
solvable
group, and $\df{\ov G_9}{\ov H_9}\simeq $ symmetric group of 9 elements .
\hfill$\qed$ for Theorem 5.0

\newpage

\subheading{\S6. The result of Veronese of order $p$}

The result for any Veronese is similar (see below), but the representation
of the proof for
$p=3$ is much more ``reader friendly''.
The result for any $p$ is as follows:

There exist 2 series $1\triangleleft A \triangleleft B \triangleleft C
\triangleleft G$ and
$ 1\triangleleft \ov A\triangleleft\ov B\triangleleft\ov C\triangleleft\ov G$
s.t.
$$ \align
   & G/C\simeq \ov G/\ov C\simeq S_{p^2}\\
 & C/B\simeq \BZ,\ \ov C/\ov B\simeq Z_q \\
& B/A\simeq \ov B/\ov A\simeq\left\{ \aligned
&(\BZ\oplus\BZ_3)^{p^2-1}\quad   p=O(3) \\
&\BZ^{p^2-1} \qquad\qquad   p\not\equiv O(3)\endaligned\right.\\
&A\simeq\ov A\simeq\left\{\aligned  \BZ_2 \quad p\ \text{odd}\\
0\quad p\ \text{even}\endaligned\right. \endalign$$\bk

\end
\Refs\widestnumber\key{MoTe12}
\ref\key A\by E.  Artin \paper Theory of braids\jour Ann. Math. \vol 48\pages
101-126\yr 1947
\endref

\ref\key CT \by J. Carlson and D. Toledo\paper Discriminant complements and
kernels of
monodromy representations\jour Duke Math. J.\toappear\endref

\ref\key DOZ\by G. Detloff, S. Orevkov and M. Zeidenberg\paper Plane curves
with a big
fundamental group of  the complement\inbook Voronezh Winter Mathematical
Schools: Dedicated
to Selim Krein, ( P. Kuchment, V. Lin, eds.), American Mathematical Society
Translations--Series 2\vol 184
  \yr 1998\endref


\ref\key FH\by Feustel / Holzapfel \paper Symmetry points and Chern
invariants of Picard-modular-surfaces\jour Math. Nach. \vol 111\pages 7-40 \yr
1983\endref


 \ref\key L \by A. Libgober   \paper Homotopy groups of the complements to
singular hypersurfaces II \jour Ann. of Math. \vol 139 \yr 1994 \pages
119-145\endref

\ref\key Mo\by B.  Moishezon  \paper
Stable branch curves and braid monodromies\inbook LNM\vol 862\pages 167-192
\endref
\ref\key MoRoTe \by B.  Moishezon, A. Robb and  M. Teicher  \paper
On  Galois covers of Hirzebruch surfaces
\jour Math. Ann. \vol305\yr1996\pages 493-539. \endref


 \ref\key MoTe1 \by B. Moishezon and M. Teicher \paper Existence of
simply connected algebraic surfaces of positive and zero indices
\jour Proceedings of the National Academy of
Sciences,  	United States of America \vol 83 \yr 1986 \pages
6665-6666
\endref
 \ref\key MoTe2 \by B. Moishezon and M. Teicher \paper Simply connected
algebraic surfaces of positive index \jour Invent. Math. \vol 89
\pages 601-643 \yr1987 \endref
\ref\key MoTe3\by B. Moishezon and M. Teicher\paper Galois coverings in
the theory of algebraic surfaces\jour  Proc. of Symp. in Pure Math.
\vol 46 \yr 1987 \pages 47-65
\endref
\ref\key MoTe4  \by B. Moishezon and M. Teicher \paper Braid group
technique in complex geometry, I, Line arrangements in $\Bbb C\Bbb P^2$ \jour
Contemp. Math. \vol 78 \yr 1988 \pages 425-555\endref

\ref\key MoTe5 \by B. Moishezon and M. Teicher \paper Finite fundamental
groups, free over ${\Bbb Z}/c{\Bbb Z}$, for Galois covers of ${\Bbb
C}\Bbb P^2$ \jour Math. Ann. \vol 293 \yr 1992 \pages 749-766
\endref


\ref\key MoTe6 \by  B. Moishezon and M. Teicher \paper Braid group
technique in complex geometry, II, From arrangements
of lines and conics to 	cuspidal curves  \inbook
Algebraic Geometry, Lecture Notes in Math. \vol 1479 	\yr 1990
\endref



\ref\key MoTe7 \by  B. Moishezon and M. Teicher \paper Braid group
techniques in complex geometry III: Projective degeneration of $V_3 $
\jour Contemp. Math.\vol 162\yr 1993\pages 313-332\endref

\ref\key MoTe8 \by B. Moishezon and M. Teicher \paper
	Braid group techniques in complex geometry IV: Braid monodromy of
the branch curve
$S_3$  of  $V_3 \rightarrow \Bbb C\Bbb P^2$ and application to
$\pi_1:(\Bbb C\Bbb P^2-S_3,*)$ \jour Contemp. Math.\vol 162 \yr 1993\pages
332-358\endref

 \ref\key MoTe9 \by B. Moishezon and M. Teicher \paper
Braid group techniques in complex geometry, V:
The fundamental group of complements of a branch curve of Veronese
generic projection \jour Communications in Analysis and Geometry \vol
4\issue 1 \yr 1996 \pages 1-120\endref

\ref\key MoTe10 \by  B. Moishezon and M. Teicher \paper
Fundamental groups of complements of curves  as
solvable groups \inbook Israel Mathematics Conference
Proceedings (AMS Publications)\vol 9\pages 329-346\yr 1996\endref



\ref\key  RoTe  \by A. Robb and  M. Teicher\paper
 Applications of braid group techniques to the decomposition
 of moduli spaces, new examples\jour Topology and its Applications\toappear
\endref

\ref\key Si\by C. Simpson  \paper A relative notion of algebraic Lie group and
applications to $n$-stacks\jour Algebraic Geometry\yr 1996\endref

\ref\key Te\by M. Teicher\paper Braid groups, algebraic surfaces and
complements of braid curves\jour PSPM 62\vol 1\yr 1997\pages 127-149 \endref


\ref\key To\by D. Toledo \paper Projective varieties with non residual finite
fundamental group\jour Extrait de Publications Mathematiques \vol 77\yr
1993\endref


\ref\key VK \by  E.R. Van Kampen \paper On the fundamental group of an
algebraic curve \jour Am. J. Math. \vol 55 \pages 255-260 \yr
1933\endref
 \ref\key Z \by O. Zariski  \book Algebraic Surfaces (Ch. VIII)
\publ Second Edition, Springer \yr 1971\endref\endRefs
